%% file: chmail.tex
\documentclass[11pt]{amsart} 
\usepackage{latexsym}                    
\usepackage{amssymb,amscd,amsthm,amsfonts,verbatim,amsopn,color}
\usepackage[dvips]{epsfig}
\newcommand{\al}{\alpha}

\newcommand{\cd}{{\mathcal D}}
\newcommand{\ci}{{\mathcal I}}

\newcommand{\Da}{\Delta}
\newcommand{\da}{\delta}

\newcommand{\hyp}{\text{Hyp}^n}
\newcommand{\hyph}{\widehat{\hyp_\lambda}}

\newcommand{\la}{\lambda}

\newcommand{\nin}{\noindent}

\newcommand{\pr}{\nin{\bf Proof. }}

\newcommand{\ra}{\rightarrow}

\newtheorem{thm}{Theorem}[section]
\newtheorem{df}  [thm]{Definition}

\newtheorem{prop}[thm]{Proposition}

\numberwithin{equation}{section}

\begin{document}

\title[Chain mails] {Directed trees in a~string, real
  polynomials with triple roots, and chain mails}

              \author{Dmitry N. Kozlov}
  \date{\today. \\[0.05cm]
 \\[0.05cm]
 \hskip15pt This research was supported by the Research Grant 
 of the Swiss Natural Science Foundation.
% \\[0.05cm]
% \hskip15pt  This paper is submitted.
}

\address{ {\it Current address:} 
Department of Mathematics, University of Bern, 
Sidler\-strasse~5, CH-3012, Bern, Switzerland.}
\address{ {\it On leave from:}
Department of Mathematics, Royal Institute of Technology, 
Stockholm, S-100 44, Sweden.  }
\address{{\it E-mail addresses:} {\tt kozlov@math-stat.unibe.ch, 
kozlov@math.kth.se.}}

\begin{abstract} 
  This paper starts with an observation that two infinite series of
  simplicial complexes, which a~priori do not seem to have anything to
  do with each other, have the same homotopy type. One series consists
  of the complexes of directed forests on a double directed string,
  while the other one consists of Shapiro-Welker models for the spaces
  of hyperbolic polynomials with a triple root.
  
  We explain this coincidence in the more general context by finding
  an explicit homotopy equivalence between complexes of directed
  forests on a double directed tree, and doubly disconnecting
  complexes of a~tree.
\end{abstract}

\maketitle

\section{Observation}

When the explicit determination of the homotopy type of families of
combinatorially defined simplicial complexes is performed, it
sometimes happens that the answers coincide. In such a~case two
natural questions arise immediately. First, one would like to have
an~explicit map between these complexes, which would provide
a~homotopy equivalence. Second, one may ask for new terminology and
a~new theorem, which would put this particular result in a~more
general context. This paper is an~example of an~analysis of
a~situation like this.

\vskip3pt

\subsection{Complexes of directed forests} $\,$
\vskip2pt

Let us start by describing the two families of simplicial complexes,
for which the computations yield the same answer.

\begin{df}\label{def}(\cite{Ko99b,Sta}).
  Let $G$ be an~arbitrary directed graph. $\Da(G)$ is the~simplicial
  complex $\Da(G)$ constructed as follows: the~vertices of $\Da(G)$
  are given by the~edges of $G$ and faces are all directed forests
  which are subgraphs of $G$.
\end{df}

The complexes of directed forests in a~given directed graph were
introduced in \cite{Ko99b} following the suggestion of R.\ Stanley. It
was proved in \cite{Ko99b} that $\Da(G)$ is shellable (hence homotopy
equivalent to a~wedge of spheres) when $G$ has a~complete sink.
Furthermore, the homotopy types of complexes $\Da(G)$ were computed
for several natural families of graphs $G$. The special case, which is
of particular importance for this paper is when $G$ is a~double
directed string.

\begin{df} (\cite{Ko99b}).
  Let $n$ be a~nonnegative integer.  A~double directed string on $n+1$
  vertices is the directed graph, denoted $L_n$, which is defined by
  $V(L_n)=[n+1]$, and $E(L_n)=\{(i\ra i+1), (i+1\ra i)\,|\,i\in[n]\}$.
\end{df}

The following proposition was also proved in \cite{Ko99b}.
    
\begin{prop} \label{prop1.3}
$$\Da(L_n)\simeq\begin{cases}
  S^{2k-1}, & \text{ if } n=3k;\\
  S^{2k}, & \text{ if } n=3k+1;\\
  \text{a point}, & \text{ if } n=3k+2.
\end{cases}$$
\end{prop}

\vskip3pt

\subsection{Spaces of monic hyperbolic polynomials with multiple roots} 
$\,$ \vskip2pt

Let us now turn to the second family of simplicial complexes. The
space of monic hyperbolic polynomials in one variable of degree~$n$,
which we denote by $\hyp$, is naturally stratified by fixing the
multiplicities of roots (in our terminology a~polynomial is called
hyperbolic if all of its roots are real). 

These strata are indexed by number partitions which refine $n$, and
for a~partition $\lambda=(\pi_1,\dots,\pi_t)$ we let $\hyph$ denote
the one-point compactification of the set of all polynomials of the
form $(x-r_1)^{\pi_1}\dots(x-r_t)^{\pi_t}$ (note that we do not
require that $r_1,\dots,r_t$ are distinct numbers).

To fix the multiplicities of roots is the same as to fix a~number
partition $\lambda\vdash n$. While one can think of a~number partition
$\lambda=(\pi_1,\dots,\pi_t)$ simply as a~set of numbers, such that
$\pi_1+\dots+\pi_t=n$, we use the notation $[\pi_1,\dots,\pi_t]$ to
denote the ordered $t$-tuple of these numbers. In such a~situation,
$[\pi_1,\dots,\pi_t]$ is called a~{\it composition} of $n$. By
forgetting the order of the numbers in the composition we get a~number
partition which is called the {\it type} of this composition. Both for
a~partition and for a~composition, we call the number of its parts the
{\it length}, and denote it by $l(\lambda)$, resp.\ 
$l([\pi_1,\dots,\pi_t])$.

The set of all compositions of $n$ is partially ordered by refinement.
Namely, let $x=[\al_1,\dots,\al_{l(x)}]$ and
$y=[\beta_1,\dots,\beta_{l(y)}]$ be two compositions of $n$, we say
that $x\leq y$ if and only if
$\al_j=\beta_{i_{j-1}+1}+\dots+\beta_{i_j}$, for $1\leq j\leq l(x)$,
and some $0=i_0<i_1<\dots<i_{l(x)}=l(y)$. Since $\beta_i>0$, for
$i=1,\dots,l(y)$, the indices $i_1,\dots,i_{l(x)-1}$ are uniquely
defined.

Given a~number partition $\la=(\pi_1,\dots,\pi_t)$ of $n$, we define
$D_\la$ to be the poset consisting of all compositions of $n$ which
are less or equal of some composition of $n$ of type $\la$.  Thus, the
number of maximal elements of $D_\la$ is equal to the number of
different ways to impose an order on the numbers $\pi_1,\dots,\pi_t$.
Note, that $D_\la$ has a~minimal element, the composition consisting
of just the number $n$, and it is easy to see that $D_\la\cup\{\hat
1\}$ is a~lattice, where $\hat 1$ is an~artificially added maximal
element.

Since the lower intervals of $D_\la$ are Boolean algebras, and $D_\la$
itself is a~meet-semilattice, there exists a~unique simplicial
complex, which we denote by $\da_\la$, such that $D_\la$ is the face
poset of $\da_\la$, i.e., the elements of $D_\la$ and the simplices of
$\da_\la$ are in bijection, and the partial order relation on $D_\la$
corresponds under this bijection to the inclusions of simplices of
$\da_\la$.

As the following theorem, proved in \cite{ShW98}, shows, the
simplicial complex $\delta_\lambda$ provides a combinatorial model for
the stratum $\hyph$.

\begin{thm} (\cite[Theorem 3.5(a)]{ShW98}).
  Let $\la$ be a~number partition of $n$, then the one-point
  compactification of the strata indexed by $\la$, $\hyph$, is
  homeomorphic to the double suspension of the simplicial complex
  $\da_\la$.
\end{thm}

Complete description of the homotopy type of the simplicial complexes
$\da_\la$ were obtained for several classes of number partitions
$\lambda$, see \cite{Ko00a,ShW98}. In particular the following is
well-known.

\begin{prop} \label{prop1.5}
(\cite[Corollary 3.10]{ShW98},\cite[Proposition 3.4(a)]{Ko00a}).
  For $\la=(k,1^t)$, where $k\geq 2$, $t\geq 0$, we have
$$\da_\la\simeq\begin{cases}
  S^{2m-1},& \text{ if } t=km, \text{ for some } m\in{\mathbb Z};\\
  S^{2m},& \text{ if } t=km+1, \text{ for some } m\in{\mathbb Z};\\
  \text{point},& \text{ otherwise}.
\end{cases}$$ 
\end{prop}

By comparing Propositions \ref{prop1.3} and \ref{prop1.5}, we get that
\begin{equation}
  \label{eq:2.1}
  \delta_{(3,1^t)}\simeq\Delta(L_t).
\end{equation}

In the next section we will provide an explicit homotopy equivalence
between $\delta_{(3,1^t)}$ and $\Delta(L_t)$, while we shall state and
prove a~more general theorem in Section~\ref{sect3}.

\section{Explicit homotopy equivalence}             

\subsection{Further descriptions of the simplicial complexes 
  $\delta_{(3,1^t)}$ and $\Delta(L_t)$} 
$\,$ \vskip2pt

Let us note an alternative description of the simplicial complex
$\delta_{(3,1^t)}$. We take as simplices all subsets $\sigma$ of
$\{1,\dots,t+2\}$, for which there exist $1\leq i\leq t+1$, such that
$i,i+1\notin\sigma$; in other words, the maximal simplices are
precisely all sets $\{1,\dots,t+2\}\setminus\{i,i+1\}$, for
$i=1,\dots,t+1$. Please note, that we step a~bit away from the usual
conventions of the simplicial complexes, in that in our description it
is not necessary that all 1-element subsets of $\{1,\dots,t+2\}$ are
simplices. For example, $\delta_{(3)}$ is simply empty, while
$\delta_{(3,1)}$ has only 2 vertices.

\vskip2pt

Also $\Delta(L_t)$ has other descriptions, two of which we list next.

\vskip3pt

\nin {\it Description 1 (via forbidden patterns).} 

\nin $\Delta(L_t)$ is the simplicial complex, with the set of vertices
being the set $\{1,\dots,2t\}$, where we take as simplices all subsets
$\sigma\subseteq\{1,\dots,2t\}$, such that
$\{i,i+1\}\not\subseteq\sigma$, for any $i=1,\dots,2t-1$. These
simplicial complexes have appeared in various guises (for example, in
\cite{Ehr} they were called {\it sparse} complexes).

\vskip3pt

With this description at hand it is obvious that the simplicial
complex $\Delta(L_t)$ is isomorphic to $\Delta(L_t')$, where $L_t'$ is
the directed graph on $t+2$ vertices, $V(L_t')=\{1,\dots,t+2\}$, given
by
$$E(L_t')=\{(1\ra 2),(t+2\ra t+1)\}\cup\{(i\ra i+1),
(i+1\ra i)\,|\,i=2,\dots,t\}.$$

$$%
  \begin{picture}(0,0)%
    \includegraphics{p1.pstex}%
  \end{picture}%
  \input{p1.pstex_t}  $$
$$\text{Figure 1. The case } t=3.$$

\vskip3pt

%\pagebreak

More generally, one has the following definition.

\begin{df}\label{df2.1}
  For an arbitrary graph $\Gamma$ we can define the {\bf independence}
  simplicial complex $\ci(\Gamma)$ as follows. The set of vertices of
  $\ci(\Gamma)$ is equal to the set of vertices of $\Gamma$. The
  subset $\sigma\subseteq\ci(\Gamma)$ is a~simplex if and only if the
  vertices in $\sigma$ can be colored with the same color.
\end{df}

Note that we use the usual requirement, that no two vertices of the
same color are connected by an edge. In the graph theory literature,
sets of vertices which can be colored with the same color are
traditionally called {\it independent sets}, which is why we chose
this name for the complexes $\ci(\Gamma)$. The condition for being
a~simplex in the Definition~\ref{df2.1} could be reformulated slightly
differently, namely, if $\sigma$ is a~simplex of $\ci(\Gamma)$, and
$(x,y)$ is an~edge of $\Gamma$, then either $x\notin\sigma$ or
$y\notin\sigma$.

\vskip3pt

\nin {\it Description 2 (via order complexes).} Recall that to any
partially ordered set (poset for short) one can associate a simplicial
complex whose vertex set is the set of the elements of the poset, and
whose set of simplices consists of all the completely ordered subsets
of elements of the poset (these ordered subsets are also known as
chains). This simplicial complex is called {\it order complex} of the
poset.  Connecting this definition to our situation, one can see that
$\Delta(L_t)$ is isomorphic to the order complex of the poset $P_t$,
given by:
\begin{itemize}
\item the set of the elements of $P_t$ is $\{1,\dots,2t\}$;
\item the partial order on $P_t$ is given by: $x>y$ if and only if
  $|x-y|\geq 2$.
\end{itemize}

$$%
  \begin{picture}(0,0)%
    \includegraphics{p2.pstex}%
  \end{picture}%
  \input{p2.pstex_t}  $$
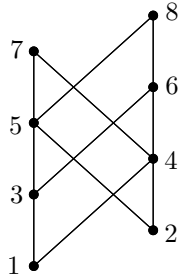
$$\text{Figure 2. The poset } P_4. $$
    
\vskip3pt
 
Note that it is rather unusual for a graph $G$ to be such that there
exists a~poset $P$, such that $\Delta(G)=\Delta(P)$. For example,
a~careful check should convince the reader that no such poset $P$
exists for the graph on the Figure~3.

$$%
  \begin{picture}(0,0)%
    \includegraphics{p25.pstex}%
  \end{picture}%
  \input{p25.pstex_t}  $$
$$\text{Figure 3.} $$

\vskip3pt

\subsection{The explicit homotopy equivalence} 
$\,$ \vskip2pt

We are now in the~position to define the map which shall give
an~explicit homotopy equivalence between $\Delta(L_t)$ and
$\delta_{(3,1^t)}$. Namely, let $\phi$ be the simplicial map defined
by
$$\begin{array}{lcccc}
\phi & : & \Delta(L_t) & \longrightarrow & \delta_{(3,1^t)}\\
&& (x\ra y)& \mapsto & x.
\end{array}$$

\begin{thm} \label{thm2.2}
  The induced map of topological spaces
  $\phi:\Delta(L_t)\longrightarrow\delta_{(3,1^t)}$ is a~homotopy
  equivalence.
\end{thm}

It is not difficult to show that $\phi$ is well-defined and prove the
Theorem~\ref{thm2.2}. However, we shall not do that here, as it will
follow from a~more general Theorem~\ref{thm3.2} in the next section.

\vskip3pt

\subsection{Chain mails} \label{ss2.3}
$\,$ \vskip2pt

To conclude this section, let us give yet another interpretation of
the simplicial complexes $\Delta(L_t)$ and $\delta_{(3,1^t)}$, which
allows to view the map $\phi$ in a~somewhat different light.

Consider $t+2$ unit rings arranged as shown on Figure~4.

$$%
  \begin{picture}(0,0)%
    \includegraphics{p3.pstex}%
  \end{picture}%
  \input{p3.pstex_t}  $$
$$\text{Figure 4.}$$

If we identify the rings with the vertices of $\delta_{(3,1^t)}$, then
the complex $\delta_{(3,1^t)}$ is identified with the simplicial
complex on the set of rings, in which the simplices are all
collections of rings with at least two missing rings in a~row.

Let us now associate a graph $\Gamma$ to this string of rings as
follows. We take vertex on each small (smaller than $\pi$) arc of the
ring, which is cut out by another ring, except for the two small arcs
on the ends. Then we connect each pair of vertices, where an edge can
be drawn without intersecting the rings; that is we connect by an edge
each pair of vertices which see each other without having to look
through the rings. Clearly we obtain a~string with $2t$ vertices, and
it is immediate that $\ci(\Gamma)=\Delta(L_t)$.

We can now see our map
$\phi:\ci(\Gamma)\longrightarrow\delta_{(3,1^t)}$ as the one which
maps each vertex of $\Gamma$ to the ring to which it belongs. This
mental image is the main reason to introduce more general structures,
which we, following the visual analogy, call chain mails.

In fact, there are various ways to generalize the situation shown on
Figure~4. We choose one which will be sufficient for our purposes.

\begin{df} \label{df2.3}
  A {\bf chain mail} is a collection of (not necessarily convex)
  polygons $R_1,\dots,R_p$ on the plane, such that
  \begin{itemize}
  \item for all $i<j$, $R_i\not\subseteq R_j$ and $R_j\not\subseteq
    R_i$;
  \item for all $i<j$, if $R_i\cap R_j$ is nonempty, then
    $\partial R_i\cap\partial R_j$ consists of two points and
    $(R_i\setminus\partial R_i)\cap(R_j\setminus\partial
    R_j)\neq\emptyset$;
  \item each triple intersection $R_i\cap R_j\cap R_k$, $i<j<k$, is
    empty.
  \end{itemize}
\end{df}

Note that if $\partial R_i\cap\partial R_j\neq\emptyset$, then the
boundary of each of the PL curves $\partial R_i$ and $\partial R_j$ is
divided into two parts: one part lies inside $R_j$, resp.\ $R_i$,
while the other one lies outside.

To each chain mail one can associate a graph $\Gamma$ as follows.
\begin{itemize}   
\item Vertices of $\Gamma$. For each pair $i\neq j$,
  $i,j\in\{1,\dots,p\}$, such that $\partial R_i\cap\partial
  R_j\neq\emptyset$, we put one vertex $v$ of the graph $G$ on the
  part of $\partial R_i$ which lies inside $R_j$ if there exists
  $k\in\{1,\dots,p\}$, $k\neq i$, $k\neq j$, such that one can connect
  $v$ with $\partial R_k$, with a~PL curve, without intersecting any
  $\partial R_h$, $h\in\{1,\dots,p\}$; see Figure~6.
\item Edges of $\Gamma$. Two vertices are connected by an~edge, if
  they can be connected by a~PL curve, which does not intersect
  $\partial R_h$, for all $h\in\{1,\dots,p\}$, other than at its
  endpoints; again see Figure~6.
\end{itemize}

\section{Homotopy equivalence in the context of doubly connected complexes}
\label{sect3}

In this section we generalize the Theorem \ref{thm2.2}, and give
a~direct proof of this generalization. 

\vskip3pt

\subsection{$k$-fold disconnecting complexes of trees}
$\,$ \vskip2pt

First we need some terminology. For any tree $T$ there is a~natural
way to augment it as follows: for each leaf of $T$ we add one extra
vertex, which is a~new leaf, and one extra edge, which connects the
new leaf with the old one. We shall call the augmented tree $\widehat
T$. Alternatively, the class of all augmented trees could be described
as the class of all trees where each path leading from a~branch vertex
to a~leaf is at least of length~2.

One can also turn the tree into a~directed tree, as one augments it,
as follows: all the edges of $T$ are replaced by pairs of edges going
in opposite directions, while the new edges are each replaced by just
one edge, which is directed to point {\it from} the new leaf. We
denote the thus obtained directed tree by $\widetilde T$.

$$%
  \begin{picture}(0,0)%
    \includegraphics{p4.pstex}%
  \end{picture}%
  \input{p4.pstex_t}  $$
$$\text{Figure 5.}$$

\vskip3pt

\begin{df} \label{df3.1}
  Given a tree $T$. For a~positive integer~$k$, let $\cd_k(T)$ be the
  simplicial complex defined as follows. The vertices of $\cd_k(T)$
  are the vertices of $\widehat T$, while $S\subseteq V(\widehat T)$
  is a~simplex if for any path $(x_1,\dots,x_t)$ connecting two
  different leaves of $\widehat T$ (note that $t\geq 4$ is
  necessitated), there exists $1\leq i\leq t-k+1$, such that
  $x_i,x_{i+1},\dots,x_{i+k-1}\notin S$.
\end{df}

We call $\cd_k(T)$ the {\bf $k$-fold disconnecting complex} of $T$.
When $k=2$ we simply call $\cd_2(T)$ the {\bf doubly disconnecting
  complex} of~$T$. Clearly, when $T$ is a~string (a~tree without
branching) with $t$ vertices, the complex $\cd_2(T)$ is isomorphic
to $\delta_{(3,1^t)}$.

It is natural to associate a chain mail to the graph $\widehat T$, as
is shown on Figure~6.

$$%
  \begin{picture}(0,0)%
    \includegraphics{p5.pstex}%
  \end{picture}%
  \input{p5.pstex_t}  $$
$$\text{Figure 6.}$$

\vskip3pt

For each vertex of $\widehat T$ we draw a (not necessarily convex)
polygon, so that these polygons form a~chain mail, i.e., satisfy
conditions of the Definition~\ref{df2.3}, and so that two polygons
intersect if and only if the corresponding vertices of $\widehat T$
are connected by an edge.  We rely on the visual clarity of the
Figure~6 for the elucidation of how such a~chain mail looks like.

However, we would like to mention one formal way to construct this
chain mail. For each vertex of $\widehat T$ take the union of the
closed halves of the edges adjacent to this vertex (in the standard
metric associated to any graph, this is the closed ball of radius 1/2
centered at this vertex). Take the closed $\epsilon$-neighborhood of
this union, where $\epsilon$ is small, and depends on the specific
imbedding of $\widehat T$ into the plane. The polygons can now be
taken as sufficiently fine (again depending on the imbedding of
$\widehat T$) PL approximations of these $\epsilon$-neighborhoods. The
graph $\Gamma$ depicted on Figure 6 is the graph associated to this
chain mail, as described in the subsection~\ref{ss2.3}. Note that
$\ci(\Gamma)=\Delta(\widetilde T)$.

\vskip3pt

\subsection{The main theorem}
$\,$
\vskip2pt

\begin{thm}\label{thm3.2}
  Let $T$ be a tree, then the simplicial complexes $\cd_2(\widehat T)$
  and $\Delta(\widetilde T)$ are homotopy equivalent.
\end{thm}

\pr The map $\phi:\Delta(\widetilde T)\longrightarrow\cd_2(\widehat
T)$ is defined analogously to the one in the previous section. One way
to describe it is to set $\phi(x\ra y)=x$. The other way, in terms of
chain mails, is to say that each vertex of the graph $\Gamma$
associated to the chain mail is mapped to the boundary of the polygon,
to which it belongs. Recall that, as mentioned above, $\ci(\Gamma)$ is
isomorphic to $\Delta(\widetilde T)$, and we are tacitly talking here
about $\phi:\ci(\Gamma)\longrightarrow\cd_2(\widetilde T)$.

\vskip2pt

\nin {\bf Claim 1.} {\it $\phi$ is well-defined.}

\vskip2pt

\nin Let $S$ be a~simplex of $\Delta(\widetilde T)$, $S\neq\emptyset$.
We need to show that $\phi(S)$ has a ``2-gap'' in each path connecting
two leaves. Without loss of generality, we can assume that $S$ is the
maximal simplex of $\Delta(\widetilde T)$. 

Assume this is not the case, and choose a path $P=(x_1,\dots,x_t)$ in
$\cd_2(\widehat T)$ which does not have a~2-gap, again $t\geq 4$. If
$(x_1\ra x_2)\notin S$, and $(x_2\ra x_3)\notin S$, then
$x_1,x_2\notin\phi(S)$, thus $P$ has a~2-gap, which is a contradiction.
In the same way, either $(x_t\ra x_{t-1})\in S$ or $(x_{t-1}\ra
x_{t-2})\in S$. So $S$ contains two edges are directed to each other,
which, by maximality of $S$, means that there must exist $2\leq i\leq
t-2$, such that $(x_{i-1}\ra x_i),(x_{i+2}\ra x_{i+1})\in S$, which in
turn implies that $x_i,x_{i+1}\notin S$. This is again a~contradiction
to the assumption that $P$ has no 2-gap.

\vskip2pt

\nin {\bf Claim 2.} {\it $\phi$ gives a homotopy equivalence.}

\vskip2pt

\nin We shall show that for each (closed) simplex
$S\subseteq\cd_2(\widehat T)$, $\phi^{-1}(S)$ is contractible. In
fact, by using induction on the number of vertices of $\widehat T$, we
shall prove that it is a cone. The theorem then follows by Quillen's
Lemma, see~\cite{Qu}.

Take any path $P=(x_1,\dots,x_t)$ in $\widehat T$, such that $P\cap
S\neq\emptyset$. Since $P\cap S$ must contain a~gap consisting of at
least two elements, we can always find, (possibly after reversing the
indexing of $x_i$'s), $2\leq i\leq t-1$, such that $x_{i-1}\in S$, and
$x_i,x_{i+1}\notin S$. 

Obviously, $(x_{i-1}\ra x_i)\in\phi^{-1}(S)$. If $\phi^{-1}(S)$ is not
a~cone with apex at $(x_{i-1}\ra x_i)$, then there must exist some
simplex $\sigma$ in $\phi^{-1}(S)$, such that $\sigma\cup\{(x_{i-1}\ra
x_i)\}\notin\phi^{-1}(S)$. This means that either $(x_i\ra
x_{i-1})\in\sigma$, or $(y\ra x_i)\in\sigma$, for some $y\neq
x_{i-1}$. The first option is ruled out by the fact that $x_i\notin
S$, furthermore, in the second option we have $y\neq x_{i+1}$, since
$x_{i+1}\notin S$. Thus, we conclude that the vertex $x_i$ has to have
valency at least~3.

Let $\widehat Q$ be the induced subtree of $\widehat T$ consisting of
all those vertices $v\in\widehat T$, which can be reached from $x_i$
without passing through $x_{i+1}$; this includes $x_i$ itself, see the
Figure~7.

$$%
  \begin{picture}(0,0)%
    \includegraphics{p6.pstex}%
  \end{picture}%
  \input{p6.pstex_t}  $$
$$\text{Figure 7.}$$

Since the valency of $x_i$ in $\widehat T$ is at least 3, its valency
in $\widehat Q$ is at least 2, hence each leaf of $\widehat Q$ is also
a~leaf of $\widehat T$. Let furthermore $\widetilde Q$ denote the
directed graph associated to $\widehat Q$; the rule of association is
the same as before: the internal edges are replaced by pairs of
oriented edges directed in opposite directions, while the edges
connecting leaves to the internal vertices of the tree are replaced by
single directed edges pointing away from the leaves. Let
$\phi_Q:\Delta(\widetilde Q)\longrightarrow\cd_2(\widehat Q)$ be the
restriction of the map~$\phi$.

By induction assumption, $\phi_Q^{-1}(S\cap\widehat Q)$ is a~cone. Let
us denote its apex by $(x\ra y)$. If $\phi^{-1}(S)$ is not a~cone with
apex $(x\ra y)$, then there must exist a~simplex
$\sigma\in\phi^{-1}(S)$, such that either $(y\ra x)\in\sigma$, or
$(z\ra y)\in\sigma$. By~the~construction of $\widehat Q$, we conclude
that the first option is impossible, and in the second option we must
have $z=x_{i+1}$. This yields a~contradiction, as it implies that
$x_{i+1}\in S$. \qed

\section{Further questions}

As mentioned above, when $T$ is a~string with $t$ vertices, $t\geq 1$,
$\cd_2(T)$ is isomorphic to $\delta_{(3,1^t)}$. More generally, one
can see that in this case, if $t\geq k-1$, $\cd_k(T)$ is isomorphic to
$\delta_{(k+1,1^{t+2-k})}$ (for the case $k=1$, $t=0$, we use the
convention that $\widehat T$ is the graph consisting of two vertices,
which are connected by an~edge).  Unfortunately, we do not know how to
generalize the Theorem \ref{thm3.2}, or even the Theorem \ref{thm2.2},
to this case.

The family of complexes, which, because of the connections to spaces
of monic hyperbolic polynomials with multiple roots, appears to us to
be of great interest, is $\{\delta_\lambda\}_{\lambda\vdash n}$. It is
effortless to extend the Definition~\ref{df3.1} to this case, so that
one can talk about $\lambda$-fold disconnecting complexes of trees.
The difficulty arises from the fact that we do not know of any natural
replacement for the complexes $\Delta(G)$.

One way to view the Theorems \ref{thm3.2} and \ref{thm2.2} is to see
them as expressing a~kind of duality. The complexes $\delta_\lambda$
are given by listing their maximal simplices, while the complexes
$\Delta(G)$ are given by listing minimal non-simplices, which turn out
to have dimension~1. Such a~description would have been impossible for
the complexes $\delta_\lambda$ themselves, since in fact they contain
complete skeletons in the dimensions up to roughly half of the number
of vertices.

In light of the above, it seems interesting to ask: {\it can one find
  simplicial complexes, having some nice combinatorial description in
  terms of minimal non-simplices of small dimension, which are
  homotopy equivalent to complexes $\delta_\lambda$?}

The Theorems \ref{thm3.2} and \ref{thm2.2} can be thought of as the
first step on the path of providing a~complete and satisfactory answer
to that question.

%1) Compute $\Delta(T)$ for various regular trees.

%2) What about complexes $\cd_k$? $k$-disconnectedness.

%3) Is there any connection to {\tt Hom}-complexes? Perhaps to $K_1^+$?

%4) What about similar models coming out of other partitions $\lambda$?

%5) What is $\phi^{-1}$, and does it have some nice description?

\end{document}

%% file: p1.pstex_t
\begin{picture}(0,0)%
\includegraphics{p1.pstex}%
\end{picture}%
\setlength{\unitlength}{3947sp}%
\begingroup\makeatletter\ifx\SetFigFont\undefined%
\gdef\SetFigFont#1#2#3#4#5{%
  \reset@font\fontsize{#1}{#2pt}%
  \fontfamily{#3}\fontseries{#4}\fontshape{#5}%
  \selectfont}%
\fi\endgroup%
\begin{picture}(4505,317)(396,285)
\put(4901,339){\makebox(0,0)[lb]{\smash{\SetFigFont{12}{14.4}{\rmdefault}{\mddefault}{\updefault}{\color[rgb]{0,0,0}$)$}%
}}}
\put(2636,339){\makebox(0,0)[lb]{\smash{\SetFigFont{12}{14.4}{\rmdefault}{\mddefault}{\updefault}{\color[rgb]{0,0,0}$\Delta($}%
}}}
\put(2207,339){\makebox(0,0)[lb]{\smash{\SetFigFont{12}{14.4}{\rmdefault}{\mddefault}{\updefault}{\color[rgb]{0,0,0}$)$}%
}}}
\put(396,339){\makebox(0,0)[lb]{\smash{\SetFigFont{12}{14.4}{\rmdefault}{\mddefault}{\updefault}{\color[rgb]{0,0,0}$\Delta($}%
}}}
\put(2371,339){\makebox(0,0)[lb]{\smash{\SetFigFont{12}{14.4}{\rmdefault}{\mddefault}{\updefault}{\color[rgb]{0,0,0}$=$}%
}}}
\put(2065,479){\makebox(0,0)[lb]{\smash{\SetFigFont{10}{12.0}{\rmdefault}{\mddefault}{\updefault}{\color[rgb]{0,0,0}4}%
}}}
\put(3877,485){\makebox(0,0)[lb]{\smash{\SetFigFont{10}{12.0}{\rmdefault}{\mddefault}{\updefault}{\color[rgb]{0,0,0}3}%
}}}
\put(721,497){\makebox(0,0)[lb]{\smash{\SetFigFont{10}{12.0}{\rmdefault}{\mddefault}{\updefault}{\color[rgb]{0,0,0}1}%
}}}
\put(2965,497){\makebox(0,0)[lb]{\smash{\SetFigFont{10}{12.0}{\rmdefault}{\mddefault}{\updefault}{\color[rgb]{0,0,0}1}%
}}}
\put(3421,491){\makebox(0,0)[lb]{\smash{\SetFigFont{10}{12.0}{\rmdefault}{\mddefault}{\updefault}{\color[rgb]{0,0,0}2}%
}}}
\put(4309,485){\makebox(0,0)[lb]{\smash{\SetFigFont{10}{12.0}{\rmdefault}{\mddefault}{\updefault}{\color[rgb]{0,0,0}4}%
}}}
\put(1171,479){\makebox(0,0)[lb]{\smash{\SetFigFont{10}{12.0}{\rmdefault}{\mddefault}{\updefault}{\color[rgb]{0,0,0}2}%
}}}
\put(4771,479){\makebox(0,0)[lb]{\smash{\SetFigFont{10}{12.0}{\rmdefault}{\mddefault}{\updefault}{\color[rgb]{0,0,0}5}%
}}}
\put(1621,479){\makebox(0,0)[lb]{\smash{\SetFigFont{10}{12.0}{\rmdefault}{\mddefault}{\updefault}{\color[rgb]{0,0,0}3}%
}}}
\end{picture}

%% file: p2.pstex_t
\begin{picture}(0,0)%
\includegraphics{p2.pstex}%
\end{picture}%
\setlength{\unitlength}{3947sp}%
\begingroup\makeatletter\ifx\SetFigFont\undefined%
\gdef\SetFigFont#1#2#3#4#5{%
  \reset@font\fontsize{#1}{#2pt}%
  \fontfamily{#3}\fontseries{#4}\fontshape{#5}%
  \selectfont}%
\fi\endgroup%
\begin{picture}(993,1704)(589,-1000)
\put(589,-1000){\makebox(0,0)[lb]{\smash{\SetFigFont{10}{12.0}{\rmdefault}{\mddefault}{\updefault}{\color[rgb]{0,0,0}1}%
}}}
\put(607,-556){\makebox(0,0)[lb]{\smash{\SetFigFont{10}{12.0}{\rmdefault}{\mddefault}{\updefault}{\color[rgb]{0,0,0}3}%
}}}
\put(601,-118){\makebox(0,0)[lb]{\smash{\SetFigFont{10}{12.0}{\rmdefault}{\mddefault}{\updefault}{\color[rgb]{0,0,0}5}%
}}}
\put(607,356){\makebox(0,0)[lb]{\smash{\SetFigFont{10}{12.0}{\rmdefault}{\mddefault}{\updefault}{\color[rgb]{0,0,0}7}%
}}}
\put(1582,599){\makebox(0,0)[lb]{\smash{\SetFigFont{10}{12.0}{\rmdefault}{\mddefault}{\updefault}{\color[rgb]{0,0,0}8}%
}}}
\put(1582,131){\makebox(0,0)[lb]{\smash{\SetFigFont{10}{12.0}{\rmdefault}{\mddefault}{\updefault}{\color[rgb]{0,0,0}6}%
}}}
\put(1576,-337){\makebox(0,0)[lb]{\smash{\SetFigFont{10}{12.0}{\rmdefault}{\mddefault}{\updefault}{\color[rgb]{0,0,0}4}%
}}}
\put(1576,-796){\makebox(0,0)[lb]{\smash{\SetFigFont{10}{12.0}{\rmdefault}{\mddefault}{\updefault}{\color[rgb]{0,0,0}2}%
}}}
\end{picture}

%% file: p25.pstex_t
\begin{picture}(0,0)%
\includegraphics{p25.pstex}%
\end{picture}%
\setlength{\unitlength}{3947sp}%
\begingroup\makeatletter\ifx\SetFigFont\undefined%
\gdef\SetFigFont#1#2#3#4#5{%
  \reset@font\fontsize{#1}{#2pt}%
  \fontfamily{#3}\fontseries{#4}\fontshape{#5}%
  \selectfont}%
\fi\endgroup%
\begin{picture}(986,686)(708,-254)
\end{picture}

%% file: p3.pstex_t
\begin{picture}(0,0)%
\includegraphics{p3.pstex}%
\end{picture}%
\setlength{\unitlength}{3947sp}%
\begingroup\makeatletter\ifx\SetFigFont\undefined%
\gdef\SetFigFont#1#2#3#4#5{%
  \reset@font\fontsize{#1}{#2pt}%
  \fontfamily{#3}\fontseries{#4}\fontshape{#5}%
  \selectfont}%
\fi\endgroup%
\begin{picture}(3922,924)(737,-671)
\put(1847,-381){\makebox(0,0)[lb]{\smash{\SetFigFont{12}{14.4}{\rmdefault}{\mddefault}{\updefault}{\color[rgb]{0,0,0}2}%
}}}
\put(2437,-376){\makebox(0,0)[lb]{\smash{\SetFigFont{12}{14.4}{\rmdefault}{\mddefault}{\updefault}{\color[rgb]{0,0,0}4}%
}}}
\put(3037,-386){\makebox(0,0)[lb]{\smash{\SetFigFont{12}{14.4}{\rmdefault}{\mddefault}{\updefault}{\color[rgb]{0,0,0}6}%
}}}
\put(3642,-401){\makebox(0,0)[lb]{\smash{\SetFigFont{12}{14.4}{\rmdefault}{\mddefault}{\updefault}{\color[rgb]{0,0,0}8}%
}}}
\put(1672,-146){\makebox(0,0)[lb]{\smash{\SetFigFont{12}{14.4}{\rmdefault}{\mddefault}{\updefault}{\color[rgb]{0,0,0}1}%
}}}
\put(2272,-136){\makebox(0,0)[lb]{\smash{\SetFigFont{12}{14.4}{\rmdefault}{\mddefault}{\updefault}{\color[rgb]{0,0,0}3}%
}}}
\put(2897,-141){\makebox(0,0)[lb]{\smash{\SetFigFont{12}{14.4}{\rmdefault}{\mddefault}{\updefault}{\color[rgb]{0,0,0}5}%
}}}
\put(3477,-136){\makebox(0,0)[lb]{\smash{\SetFigFont{12}{14.4}{\rmdefault}{\mddefault}{\updefault}{\color[rgb]{0,0,0}7}%
}}}
\end{picture}

%% file: p4.pstex_t
\begin{picture}(0,0)%
\includegraphics{p4.pstex}%
\end{picture}%
\setlength{\unitlength}{3947sp}%
\begingroup\makeatletter\ifx\SetFigFont\undefined%
\gdef\SetFigFont#1#2#3#4#5{%
  \reset@font\fontsize{#1}{#2pt}%
  \fontfamily{#3}\fontseries{#4}\fontshape{#5}%
  \selectfont}%
\fi\endgroup%
\begin{picture}(6673,2375)(941,-1936)
\put(4051,-1936){\makebox(0,0)[lb]{\smash{\SetFigFont{12}{14.4}{\rmdefault}{\mddefault}{\updefault}{\color[rgb]{0,0,0}$T$}%
}}}
\put(6376,-1936){\makebox(0,0)[lb]{\smash{\SetFigFont{12}{14.4}{\rmdefault}{\mddefault}{\updefault}{\color[rgb]{0,0,0}$\widetilde T$}%
}}}
\put(1501,-1936){\makebox(0,0)[lb]{\smash{\SetFigFont{12}{14.4}{\rmdefault}{\mddefault}{\updefault}{\color[rgb]{0,0,0}$\widehat T$}%
}}}
\end{picture}

%% file: p5.pstex_t
\begin{picture}(0,0)%
\includegraphics{p5.pstex}%
\end{picture}%
\setlength{\unitlength}{3947sp}%
\begingroup\makeatletter\ifx\SetFigFont\undefined%
\gdef\SetFigFont#1#2#3#4#5{%
  \reset@font\fontsize{#1}{#2pt}%
  \fontfamily{#3}\fontseries{#4}\fontshape{#5}%
  \selectfont}%
\fi\endgroup%
\begin{picture}(2419,2577)(744,-1977)
\end{picture}

%% file: p6.pstex_t
\begin{picture}(0,0)%
\includegraphics{p6.pstex}%
\end{picture}%
\setlength{\unitlength}{3947sp}%
\begingroup\makeatletter\ifx\SetFigFont\undefined%
\gdef\SetFigFont#1#2#3#4#5{%
  \reset@font\fontsize{#1}{#2pt}%
  \fontfamily{#3}\fontseries{#4}\fontshape{#5}%
  \selectfont}%
\fi\endgroup%
\begin{picture}(2734,1855)(44,-1178)
\put(1399,-205){\makebox(0,0)[lb]{\smash{\SetFigFont{10}{12.0}{\rmdefault}{\mddefault}{\updefault}{\color[rgb]{0,0,0}$x_i$}%
}}}
\put(2425,-931){\makebox(0,0)[lb]{\smash{\SetFigFont{10}{12.0}{\rmdefault}{\mddefault}{\updefault}{\color[rgb]{0,0,0}$\widehat Q$}%
}}}
\put(757, 29){\makebox(0,0)[lb]{\smash{\SetFigFont{10}{12.0}{\rmdefault}{\mddefault}{\updefault}{\color[rgb]{0,0,0}$x_{i-1}$}%
}}}
\put(2017,-223){\makebox(0,0)[lb]{\smash{\SetFigFont{10}{12.0}{\rmdefault}{\mddefault}{\updefault}{\color[rgb]{0,0,0}$x_{i+1}$}%
}}}
\end{picture}